\documentclass{svmult}
\usepackage[bottom]{footmisc}% places footnotes at page bottom

\usepackage{mathptmx}
\usepackage{amsmath,amssymb}
\usepackage{mathrsfs}
\usepackage{microtype}
\usepackage{upref}
\usepackage{url}
\let\mathcal\mathscr
\usepackage[mathscr]{eucal}
\usepackage[all]{xy}
\SelectTips{cm}{}

\smartqed

\allowdisplaybreaks[4]

\newcommand{\cprime}{\/{\mathsurround=0pt$'$}}

\newcommand*{\pd}[2]{\mathchoice{\frac{\partial#1}{\partial#2}}
  {\partial#1/\partial#2}{\partial#1/\partial#2}
  {\partial#1/\partial#2}}
\newcommand*{\fd}[2]{\mathchoice{\frac{\delta#1}{\delta#2}}
  {\delta #1/\delta#2}{\delta#1/\delta#2}{\delta#1/\delta#2}}

\newcommand{\envert}[2][\right]{\relax
  \ifx#1\right\relax \left\lvert\else#1\lvert\fi#2#1\rvert}
\let\abs=\envert

\let\phi\varphi

\newcommand{\id}{\mathrm{id}}

\DeclareMathOperator{\CL}{CL}

\newcommand{\Ev}{E}

\providecommand{\href}[2]{#2}

\newcommand*{\email}[1]{\href{mailto:#1}{\begingroup \urlstyle{rm}\Url{#1}}}

\makeatletter
\renewcommand{\@fnsymbol}[1]{}
\makeatother

\begin{document}

\title*{Hamiltonian structures for general PDEs\thanks{This work was supported
    in part by the NWO-RFBR grant 047.017.015 (PK, IK, and AV),
    RFBR-Consortium E.I.N.S.T.E.IN grant 06-01-92060 (IK, AV, and RV) and
    RFBR-CNRS grant 08-07-92496 (IK and AV).}}
  
  \author{P. Kersten\inst{1}
  \and I.S. Krasil\cprime shchik\inst{2} \and A.M. Verbovetsky\inst{2} \and R.
  Vitolo\inst{4} } \institute{P.H.M. Kersten, University of Twente, Postbus
  217, 7500 AE
  Enschede, the Netherlands\\
  \email{kerstenphm@ewi.utwente.nl} \and I.S. Krasil\cprime shchik, A.M.
  Verbovetsky, Independent University of Moscow, B. Vlasevsky 11, 119002
  Moscow, Russia\\
  \email{josephk@diffiety.ac.ru, verbovet@mccme.ru} \and R. Vitolo, Dept.\ of
  Mathematics ``E. De Giorgi'', Universit\`a del
  Salento, via per Arnesano, 73100 Lecce, Italy\\
  \email{raffaele.vitolo@unile.it} } \titlerunning{Hamiltonian structures for
  general PDEs} \authorrunning{P.H.M. Kersten, I.S. Krasil\cprime shchik, A.M.
  Verbovetsky, R. Vitolo}
  \maketitle

\begin{abstract}
  We sketch out a new geometric framework to construct Hamiltonian operators
  for generic, non-evolutionary partial differential equations. Examples on
  how the formalism works are provided for the KdV equation, Camassa-Holm
  equation, and Kupershmidt's deformation of a bi-Hamiltonian system.
\end{abstract}

\section{Introduction}

In this short paper we will discuss the following question: What happens to a 
Hamiltonian operator of an evolution system if we change coordinates so that 
the system becomes non-evolution?  

Using the traditional definition of a Hamiltonian structure one cannot answer
this question, since the definition is tied to evolution form of the system at
hand.  However, first, not all equations have a natural evolution form, and,
second, an evolution form of a system of equations is not unique.  Let us
consider some examples.

\begin{example}[KdV]
  It is well known that the KdV equation $u_t=u_{xxx}+6uu_x$ has two compatible
Hamiltonian operators:
\begin{equation}
  \label{eq:2}
  A_1=D_x,\qquad A_2=D_{xxx}+4uD_x+2u_x,
\end{equation}
so that the equation can be written in the following ways:
\begin{align*}
  u_t=u_{xxx}+6uu_x&=D_x\fd{}{u}(u^3-u_x^2/2)\\
  &=(D_{xxx}+4uD_x+2u_x)\fd{}{u}(u^2/2),
\end{align*}
where $\fd{}{u}$ denotes the Euler operator (the variational derivative) and
is applied to the two Hamiltonian densities.

Let us introduce new dependent variables $v$ and $w$ and rewrite the KdV
equation in the form
\begin{equation}
  \label{eq:2x}
  u_x=v, \quad
  v_x=w, \quad
  w_x=u_t-6uv.
\end{equation}
In the new coordinates, the KdV still has an evolutionary form, but with
respect to another independent variable ($x$ instead of~$t$). A natural
question arises then: Is the KdV equation in the form~\eqref{eq:2x}
Hamiltonian? An affirmative answer to this question was obtained by Tsarev
in~\cite{Tsarev:HPSIEqCMMMP}.  He proved that transformations of the
type~\eqref{eq:2x} preserve the Hamiltonian property of all evolution systems
for which the Cauchy problem is solvable.  Our approach is very different from
Tsarev's one.  Below we explain why this fact holds true for all
transformations of variables and without the assumption on the Cauchy problem.
We will also show how to compute the Hamiltonian structure in new coordinates.
For the above example the answer is the following:
\begin{equation}
  \label{a:99}
  \begin{aligned}
    \begin{pmatrix}
      u \\
      v \\
      w
    \end{pmatrix}_x&=\begin{pmatrix}
      0 & -1 & 0  \\
      1 & 0  & -6u \\
      0 & 6u & D_t
    \end{pmatrix}\begin{pmatrix}
      \fd{}{u} \\
      \fd{}{v} \\
      \fd{}{w}
    \end{pmatrix}(uw-v^2/2+2u^3)
    \\
    &=\begin{pmatrix}
      0 & -2u & -D_t-2v  \\
      2u & D_t  & -12u^2-2w \\
      -D_t+2v & 12u^2+2w & 8uD_t+4u_t
    \end{pmatrix}\begin{pmatrix}
      \fd{}{u} \\
      \fd{}{v} \\
      \fd{}{w}
    \end{pmatrix}(-3u^2/2-w/2).
  \end{aligned}
\end{equation}
\end{example}

\begin{example}[Camassa-Holm equation]
  \label{e:1}
  Camassa and Holm have written their equation 
  $u_t-u_{txx}-uu_{xxx}-2u_xu_{xx}+3uu_x=0$ in a bi-Hamiltonian form by 
  introducing the new variable $m=u-u_{xx}$.  The equation now takes the form
  \begin{equation}\label{eq:8}
    m_t=-um_x-2u_xm=B_1\,\fd{\mathcal{H}_1}{m}=B_2\,\fd{\mathcal{H}_2}{m}
  \end{equation}
  with
  \begin{gather*}
    B_1=-(mD_x+D_xm),\quad \mathcal{H}_1=\frac12\int mu\,dx, \\
    B_2=D_x^3-D_x,\quad\mathcal{H}_2=\frac12\int(u^3+uu_x^2)\,dx.
  \end{gather*}
  Note that $\mathcal{H}_1$ and $\mathcal{H}_2$ are viewed as functionals in
  $m$ and $u$, but not in $u$ solely. To get rid of $m$, one is forced to
  assume that $u=(1-D_x^2)^{-1}m$ in the Hamiltonian densities. The use of the
  inverse of the operator $1-D_x^2$ is not elegant from mathematical
  viewpoint.  We will find a bi-Hamiltonian structure for the Camassa-Holm
  equation written in the initial non-evolution form and thus get rid of the
  term $(1-D_x^2)^{-1}$.
\end{example}

\begin{example}[Kupershmidt deformation]
  Consider a bi-Hamiltonian evolution system of equations
  $u_t=f(t,x,u,u_x,u_{xx},\dots)$, $u$ and $f$ being vector functions, with
  compatible Hamiltonian operators~$A_1$ and~$A_2$ and a Magri hierarchy of
  conserved densities $H_1$,~$H_2$,~\dots
\begin{equation*}
  D_t(H_i)=0,\qquad A_1\fd{H_i}{u}=A_2\fd{H_{i+1}}{u}.
\end{equation*}
In~\cite{Kupershmidt:KdV6}, Kupershmidt defined what he called the
\emph{nonholonomic deformation} of the above system:
\begin{equation}
  \label{eq:1}
  u_t=f-A_1(w),\qquad A_2(w)=0.
\end{equation}
We call system~\eqref{eq:1} the \emph{Kupershmidt deformation} of the system
$u_t=f$. The motivating example of this construction is the so-called KdV6
equation (see~\cite{Karasu(Kalkani)KarasuSakovichSakovichTurhan:ANIGKdV})
\begin{equation}
  \label{eq:5}
  u_t=u_{xxx}+6uu_x-w_x,\quad w_{xxx}+4uw_x+2u_xw=0
\end{equation}
which is the Kupershmidt deformation of the KdV equation.  The authors
of~\cite{Karasu(Kalkani)KarasuSakovichSakovichTurhan:ANIGKdV} have shown that
the KdV6 passes the Painlev\'e test and conjectured that the system is
integrable.  Kupershmidt, in~\cite{Kupershmidt:KdV6}, found a hierarchy of
conservation laws of the KdV6 as a particular case of the following general
fact.  {\def\thetheorem{}
\begin{theorem}[Kupershmidt]
  Let $u_t=f$ be an evolution bi-Hamiltonian system\textup{,} with
  $A_1$\textup{,}~$A_2$ being the corresponding Hamiltonian operators.  If
  this equation has a Magri hierarchy of conserved densities
  $\frac{dH_i}{dt}=0$\textup{,} $A_1\frac{\delta H_i}{\delta u}
  =A_2\frac{\delta H_{i+1}}{\delta u}$ then $H_1$\textup{,}~$H_2$,~\dots\ are
  conserved densities for~\eqref{eq:1}.
\end{theorem}}
\addtocounter{theorem}{-1}
\begin{proof}
  \begin{multline*}
    \frac{dH_i}{dt} =\left\langle\frac{\delta H_i}{\delta
        u},f+A_1(w)\right\rangle =\left\langle -A_1\frac{\delta H_i}{\delta
        u},
      w\right\rangle \\
    =\left\langle -A_2\frac{\delta H_{i+1}}{\delta u}, w\right\rangle
    =\left\langle\frac{\delta H_{i+1}}{\delta
        u},A_2(w)\right\rangle=0.\quad\qed
  \end{multline*}
\end{proof}
Kupershmidt also conjectured that $H_1,$~$H_2$,~\ldots\ commute in some sense
so that the KdV6 is indeed integrable.  Below we will see that this is true
and, moreover, system~\eqref{eq:1} is bi-Hamiltonian.
\end{example}

Our framework to study Hamiltonian structures for general PDEs is the geometry
of jet spaces and differential equations.  We assume the reader to be familiar
with the geometric approach to differential equations and hence we include
only the notation and the coordinate descriptions in the next section.  We
refer the reader to the
books~\cite{KrasilshchikVinogradov:SCLDEqMP,KrasilshchikVerbovetsky:HMEqMP}
for further information.

\section{Notation: infinite jets and differential equations}
In what follows everything is supposed to be smooth.

We denote an infinite jet space by~$J^\infty$.  This can be the space of jets
of submanifolds, maps, sections of a bundle, and so on, and it is not
important to us here.  Coordinates on~$J^\infty$ are $x_i$ (independent
variables, $i=1,\ldots,n$) and $u^j_{\sigma}$ (dependent variables,
$j=1,\ldots,m$, $\sigma$ being multi-indices).

The formulas
\begin{equation*}
  D_i=\pd{}{x_i} + \sum_{j,\sigma} u_{\sigma i}^j\pd{}{u_{\sigma}^j}
\end{equation*}
provide expressions for the total derivatives.  The vector fields~$D_i$ span
the Cartan distribution on~$J^\infty$.  To every vector function
on~$J^\infty$, there corresponds the evolutionary field
\begin{equation*}
  \Ev_\phi=\sum_{j,\sigma}D_\sigma(\phi^j)\pd{}{u_\sigma^j}.
\end{equation*}

The matrix differential operator
\begin{equation*}
  \ell_f = \left\Vert\sum_{\sigma}\pd{f^i}{u^j_\sigma}D_\sigma\right\Vert.
\end{equation*}
is the linearization of a vector function~$f$. It is defined by the formula
$\ell_f(\phi)=\Ev_\phi(f)$.  The linearization is a differential
operator in total derivatives; we shall call such operators
\emph{$\mathcal{C}$-differential operators}.

The coordinate expression for the adjoint $\mathcal{C}$-differential operator
is
\begin{equation*}
  \Delta^*=\left\Vert\sum_{\sigma}(-1)^{\abs{\sigma}}
    D_\sigma a^{ji}_\sigma\right\Vert
\end{equation*}
if $\Delta=\left\Vert\sum_{\sigma}a^{ij}_\sigma D_{\sigma}\right\Vert$.

Let $F_k(x_i,u_{\sigma}^j)=0$, $k=1$,~\ldots~$l$, be a system of differential
equations.  Then the relations $F=(F_1,\dots,F_l)=0$ together with
$D_{\sigma}(F)=0$ define its infinite prolongation~$\mathcal{E}\subset
J^{\infty}$.  For the sake of brevity we shall call the infinite prolongation
of a system of differential equations the equation. The operator
$\ell_{\mathcal{E}}=\ell_F|_{\mathcal{E}}$ is the linearization of the
equation~$\mathcal{E}$.

In this paper, we only consider equations~$\mathcal{E}$ whose
linearization~$\ell_{\mathcal{E}}$ is 
\emph{normal} in the following sense.
\begin{definition}
A $\mathcal{C}$-differential operator~$\nabla$ called 
\emph{normal} if  the compatibility operators for both $\nabla$ and
$\nabla^*$ are trivial.  In other words, if there exists a
$\mathcal{C}$-differential operator~$\Delta$ such that
$\Delta\circ\nabla=0$ on~$\mathcal{E}$ then $\Delta=0$
on~$\mathcal{E}$ as well, and the same holds true with $\nabla^*$
instead of $\nabla$.
\end{definition}

An evolutionary field $\Ev_{\phi}$ is a symmetry of the equation~$\mathcal{E}$
if $\Ev_{\phi}(F)|_{\mathcal{E}}=\ell_{\mathcal{E}}(\phi)=0$. If
$\Ev_{\phi}$ is a symmetry then $\phi$ is said to be its generating
function.  We often identify symmetries with their generating functions.

A vector function $S=(S^1,\dots,S^n)$ on~$\mathcal{E}$ is a conserved current 
if $\sum_iD_i(S^i)=0$ on~$\mathcal{E}$.  A conserved current is trivial if 
there exist functions~$T_{ij}$ on~$\mathcal{E}$ such that 
$S^i=\sum_{j<i}D_j(T^{ji})-\sum_{i<j}D_j(T^{ij})$.  
  
Conservation laws of~$\mathcal{E}$ are classes of conserved currents modulo
trivial ones.  To every conservation law, there correspond its generating
function, which is computed in the following way.  If $S=(S^1,\dots,S^n)$ is a
conserved current, so that $\sum_iD_i(S^i)=0$~on~$\mathcal{E}$, then there
exists a $\mathcal{C}$-differential operator~$\Delta$ such that
$\sum_iD_i(S^i)=\Delta(F)$~on~$J^\infty$.  The generating function of the
conservation law is defined by $\psi=(\psi_1,\dots,\psi_m)=\Delta^*(1)$.  Note
that $\psi=0$ if and only if the conserved current~$S$ is trivial.  One can
prove that every generating function~$\psi$ satisfies the equation
$\ell^*_{\mathcal{E}}(\psi)=0$, so that the set $\CL(\mathcal{E})$ of
conservation laws of~$\mathcal{E}$ is a subset in the kernel
of~$\ell_{\mathcal{E}}^*$, $\CL(\mathcal{E})\subset\ker\ell_{\mathcal{E}}^*$.

\section{Cotangent bundle to an equation}

Let us introduce our main hero.  For every differential equation~$\mathcal{E}$ 
we define a canonical covering 
$\tau^*\colon\mathcal{L}^*(\mathcal{E})\to\mathcal{E}$, called the
\emph{$\ell^*$-covering}.  The equation~$\mathcal{L}^*(\mathcal{E})$ 
is given by the system
\begin{equation*}
  \ell_F^*(p)=0,\quad F=0,
\end{equation*}
if $\mathcal{E}$ is given by $F=0$.  Here $p=(p^1,\ldots,p^l)$ are new
dependent variables, $l$ being the number of equations $F=(F_1,\dots,F_l)$.
We endow $\mathcal{L}^*(\mathcal{E})$ with the structure of a supermanifold by
choosing the variables~$p^k$ to be odd.  The covering $\tau^*$ is the natural
projection $\tau^*\colon(u^j_{\sigma},p^k_{\sigma})\mapsto(u^j_{\sigma})$.

Note that
\begin{equation}
  \label{a:2}
  \langle F,p\rangle=\sum_{i=1}^{l}F_ip^i
\end{equation}
is the Lagrangian for the equation~$\mathcal{L}^*(\mathcal{E})$.

It is easily shown that $\ell_{\mathcal{L}^*(\mathcal{E})}$ is normal if 
$\ell_{\mathcal{E}}$ is normal.

From the above definition it is not seen why we said that $\ell^*$-covering is
\emph{canonical}.  Indeed, the definition uses the embedding $\mathcal{E}\to
J^\infty$, but later we will show that $\mathcal{L}^*(\mathcal{E})$ is
independent of the choice of this embedding.

\begin{remark}
  \label{a:5}
  For an arbitrary $\mathcal{C}$-differential operator~$\Delta$ one can define
  the $\Delta$-covering in the same way as the $\ell^*$-covering is associated
  with the operator $\ell_{\mathcal{E}}^*$.
\end{remark}

The most interesting for us property of the $\ell^*$-covering is given by the 
following theorem.

\begin{theorem}
  \label{a:1}
  There is a natural $1$-$1$ correspondence between the symmetries
  of~$\mathcal{E}$ and the conservation laws of~$\mathcal{L}^*(\mathcal{E})$ 
  linear along the fibers of~$\tau^*$.
\end{theorem}
The expression ``linear conservation law'' means that the corresponding
conserved current is linear along the fibers of~$\tau^*$ (i.e., linear in
variables~$p^k$).  Here and below we skip the proofs that can be found
in our joint paper with S. Igonin
\cite{IgoninKerstenKrasilshchikVerbovetskyVitolo:VBGPDE}.  Let us
nevertheless describe the correspondence stated in the theorem in terms of
generating functions.  If $\phi$ is a symmetry of equation~$\mathcal{E}$ then
there exists a $\mathcal{C}$-differential operator~$\Delta$ such that
$\ell_F(\phi)=\Delta(F)$. Consider the adjoint operator~$\Delta^*$.  It can be
naturally identified with a fiberwise linear vector function~$\phi_{\Delta}$
on~$\mathcal{L}^*(\mathcal{E})$.  Then the vector function
$(\phi,\phi_{\Delta})$ is the generating function of the conservation
law that corresponds to the symmetry~$\phi$.

In the geometry of differential equation it is very useful to construct an
analogy with geometry of finite dimensional manifolds.  We shall now use this
approach to clarify the meaning of the above theorem. Let us start building
our analogy with the following two rather standard correspondences
(cf.~\cite{Vinogradov:CAnPDEqSC} and references therein):
\begin{equation*}
  \begin{aligned}
    \text{\textbf{Manifold~$M$}}\quad\phantom{\longleftrightarrow}
    \quad&\text{\textbf{PDE~$\mathcal{E}$}} \\[1ex]
    \text{functions}\quad\longleftrightarrow\quad&\text{conservation laws} \\
    \text{vector fields}\quad\longleftrightarrow\quad&\text{symmetries}
  \end{aligned}
\end{equation*}
Now, using Theorem~\ref{a:1}, we can say that the analog of the
$\ell^*$-covering is a vector bundle such that vector fields on the base are
in $1$-$1$ correspondence with fiberwise linear functions on the total space
of the bundle.  Obviously, such a bundle is the cotangent bundle.  So, the
$\ell^*$-covering is the cotangent bundle to an equation, and we can continue
our manifold-equation dictionary:
\begin{equation*}
  \begin{aligned}
    \text{\textbf{Manifold~$M$}}\quad\phantom{\longleftrightarrow}
    \quad&\text{\textbf{PDE~$\mathcal{E}$}} \\[1ex]
    \text{functions}\quad\longleftrightarrow\quad&\text{conservation laws} \\
    \text{vector fields}\quad\longleftrightarrow\quad&\text{symmetries} \\
    \text{$T^*(M)$}\quad\longleftrightarrow\quad&\text{$\mathcal{L}^*
    (\mathcal{E})$}
    \end{aligned}
  \end{equation*}

\begin{remark}
  This dictionary can be easily extended:
  \begin{equation*}
    \begin{aligned}
      \text{\textbf{Manifold~$M$}}\quad\phantom{\longleftrightarrow}
      \quad&\text{\textbf{PDE~$\mathcal{E}$}} \\[1ex]
      \text{functions}\quad\longleftrightarrow\quad&\text{conservation laws} \\
      \text{vector fields}\quad\longleftrightarrow\quad&\text{symmetries} \\
      \text{$T^*(M)$}\quad\longleftrightarrow\quad&\text{$\mathcal{L}^*
        (\mathcal{E})$} \\
      \text{$T(M)$}\quad\longleftrightarrow\quad&\text{$\mathcal{L}
        (\mathcal{E})$} \\
      \text{De Rham complex}\quad\longleftrightarrow\quad&E_1^{0,n-1}\to
      E_1^{1,n-1}\to E_1^{2,n-1}\to\cdots
    \end{aligned}
  \end{equation*}
  Here $\mathcal{L}(\mathcal{E})$ is the $\ell$-covering
  (see Remark~\ref{a:5}).  The complex 
  $E_1^{0,n-1}\to E_1^{1,n-1}\to E_1^{2,n-1}\to\cdots$ is $(n-1)$st line of 
  the Vinogradov $\mathcal{C}$-spectral sequence 
  (see~\cite{Vinogradov:CAnPDEqSC} and references therein).  In this paper we 
  use only the first three entries of the dictionary.
\end{remark}

\begin{remark}
  In~\cite{Kupershmidt:GJBSLHF}, Kupershmidt defined the cotangent bundle to a
  bundle.  This construction can be identified with the $\ell^*$-covering of
  the system
  \begin{equation*}
    u^1_t=0,\quad u^2_t=0,\quad\dots\quad u^m_t=0.
  \end{equation*}
\end{remark}

At this point, a natural question may arise: what is the analog of the
Poisson bracket on the cotangent bundle?  The answer is that the
$\ell^*$-covering is endowed with a canonical Poisson bracket.  More
precisely, since we changed the parity of fibers in the $\ell^*$-covering,
this bracket is a superbracket and is the analog of the Schouten bracket.
We shall call it the \emph{variational Schouten bracket}.

To define the bracket, recall that $\mathcal{L}^*(\mathcal{E})$ has the
Lagrangian structure~\eqref{a:2}.  Hence, by the Noether theorem there is a
$1$-$1$ correspondence between conservation laws
on~$\mathcal{L}^*(\mathcal{E})$ and Noether symmetries
of~$\mathcal{L}^*(\mathcal{E})$.  If $\psi$ is the generating function of a
conservation law, then $\Ev_{\psi}$ is the corresponding Noether symmetry.
The set of Noether symmetries is a Lie superalgebra with respect to the
commutator, so we obtain a structure of Lie superalgebra on conservation laws
on~$\mathcal{L}^*(\mathcal{E})$ uniquely determined by the equality
\begin{equation}
  \Ev_{[\![\psi_1,\psi_2]\!]}=[\Ev_{\psi_1},\Ev_{\psi_2}].
\end{equation}

According to our manifold-equation dictionary, conservation laws
on~$\mathcal{L}^*(\mathcal{E})$ correspond to functions on~$T^*(M)$.  The
latter are skew multivectors on~$M$ (this is why we have changed the parity of
fibers of the $\ell^*$-covering---to get skew-symmetric multivectors).  So, we
shall call conservation laws on~$\mathcal{L}^*(\mathcal{E})$ the
\emph{variational multivectors}.  Linear conservation laws, as we saw, are
vectors, biliner ones are bivectors and so~on.

The generating function of a variational $k$-vector is a vector function on
$\mathcal{L}^*(\mathcal{E})$ which is $(k-1)$-linear along $\tau^*$-fibers.
Such a function can be identified with a $(k-1)$-linear
$\mathcal{C}$-differential operator on~$\mathcal{E}$.  In coordinates, this
correspondence boils down to the change $p_{\sigma}\mapsto D_{\sigma}$.  Thus,
we can (and will) identify variational multivectors to multilinear
$\mathcal{C}$-differential operators.

More precisely, in the above identification we will use not operators but
equivalence classes of $\mathcal{C}$-differential operators modulo operators
divisible by~$\ell^*_{\mathcal{E}}$.  This is being done, because operators of
the form $\square\circ\ell^*_{\mathcal{E}}$ correspond to trivial functions
on~$\mathcal{L}^*(\mathcal{E})$.  But we will not change terminology, we say
operator instead of the equivalence class.

For the sake of brevity and because we are interested in the Hamiltonian
formalism, let us restrict ourselves to bivectors, which are identified with
linear $\mathcal{C}$-differential operators.  Formulas presented below for
bivectors ($=$~linear operators) can be easily generalised to multivectors
($=$~multilinear operators).

\begin{theorem}
  An operator~$A$ is a variational bivector on equation~$\mathcal{E}$ if and
  only if it satisfies the condition
  \begin{equation*}
      \ell_{\mathcal{E}} A=A^*\ell^*_{\mathcal{E}}.
    \end{equation*}
\end{theorem}

\begin{remark}
  If $\mathcal{E}$ is written in evolution form then the above condition 
  implies that $A^*=-A$.
\end{remark}

From this theorem it follows that a Hamiltonian operator~$A$ takes a 
generating function of a conservation law~$\psi$ to a symmetry~$A(\psi)$.

This is the formula for the variational Schouten bracket of two bivectors:
\begin{equation*}
  \begin{aligned}
    &[\![A_1,A_2]\!](\psi_1,\psi_2) \\[1ex]
    &\qquad=\ell_{A_1,\psi_1}(A_2(\psi_2))
    -\ell_{A_1,\psi_2}(A_2(\psi_1)) \\[1ex]
    &\qquad+\ell_{A_2,\psi_1}(A_1(\psi_2))
    -\ell_{A_2,\psi_2}(A_1(\psi_1)) \\[1ex]
    &\qquad\qquad-A_1(B_2^*(\psi_1,\psi_2))-A_2(B_1^*(\psi_1,\psi_2)),
  \end{aligned}
\end{equation*}
where $\ell_{A,\psi}=\ell_{A(\psi)}-A\ell_{\psi}$ and the operators $B_i^*$
are defined by the equalities:
\begin{equation*}
  \begin{aligned}
    \ell_FA_i-A_i^*\ell_F^*&=B_i(F,\cdot)\quad\text{on $J^{\infty}$}, \\
B_i^*(\psi_1,\psi_2)&=B_i^{*_1}(\psi_1,\psi_2)|_{\mathcal{E}}.
\end{aligned}
\end{equation*}
Here ${}^{*_1}$ denotes that the adjoint operator is computed with respect to 
the first argument.  The operators $B_i^*$ are skew-symmetric and skew-adjoint 
in each argument.  Note that if $\mathcal{E}$ is in evolution form then
$B_i^*(\psi_1,\psi_2)=\ell^*_{A_i,\psi_2}(\psi_1)$.  

Now we are in position to give a definition of a Hamiltonian structure for a 
general PDE.

\begin{definition}
  A variational bivector~$A$ is called \emph{Hamiltonian} if $[\![A,A]\!]=0$.
\end{definition}

A Hamiltonian bivector~$A$ gives rise to a Poisson bracket
\begin{equation}
  \{\psi_1,\psi_2\}_A=\Ev_{A(\psi_1)}(\psi_2)+\Delta^*(\psi_2),
\end{equation}
where $\psi_1$ and $\psi_2$ are conservation laws of~$\mathcal{E}$ and the 
operator~$\Delta$ is defined by the relation
$\ell_F(A(\psi_1))=\Delta(F)$.

As in the evolution case, we call an equation \emph{bi-Hamiltonian} if it
possesses two Hamiltonian structures $A_1$ and $A_2$ such
that~$[\![A_1,A_2]\!]=0$.

An infinite series of conservation laws $\psi_1$, $\psi_2$,~\ldots is called a 
\emph{Magri hierarchy} if for all $i$ we have $A_1(\psi_i)=A_2(\psi_{i+1})$.  
In the standard way one can show that
$\{\psi_i,\psi_j\}_{A_1}=\{\psi_i,\psi_j\}_{A_2}=0$ for all $i$ and $j$.  

Now let us return to the question of invariance of the $\ell^*$-covering.
Suppose the equation~$\mathcal{E}$ under consideration is embedded in two
different jet spaces
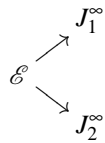
\begin{equation*}
   \xymatrixcolsep{1pc}
    \xymatrixrowsep{0.5pc}
    \xymatrix{
      &J_1^{\infty} \\
      \mathcal{E}\ar[ur]\ar[dr] \\
      &J_2^{\infty}
    }
\end{equation*}
We encountered an example of this situation when discussed the KdV equation,
with $J_1^{\infty}$ being jets with coordinates $x$, $t$ and $u$, while
$J_2^{\infty}$ being jets with coordinates $x$, $t$, $u$, $v$, and~$w$.  Now,
we have two linearization operators, $\ell_{\mathcal{E}}^1$ and
$\ell_{\mathcal{E}}^2$, the former computed using the embedding
$\mathcal{E}\to J^{\infty}_1$ and the latter is obtained using the embedding
$\mathcal{E}\to J^{\infty}_2$.  It is not difficult to show that these two
linearization operators are related by the following diagram:
\begin{equation}
  \label{a:3}
  \xymatrixcolsep{5pc}
  \xymatrix{
    \bullet\ar[r]_{\ell_{\mathcal{E}}^1}\ar@<.5ex>[d]^{\alpha}
    &\bullet\ar@<.5ex>[d]^{\alpha'}\ar@/_1pc/@<-1ex>[l]_{s_1} \\
    \bullet\ar@<.5ex>[u]^{\beta}\ar[r]^{\ell_{\mathcal{E}}^2}
    &\bullet\ar@<.5ex>[u]^{\beta'}\ar@/^1pc/@<1ex>[l]^{s_2}
    }
\end{equation}
where all arrows are $\mathcal{C}$-differential operators on~$\mathcal{E}$ 
satisfying the following relations:
\begin{equation}
  \label{a:4}
  \ell_{\mathcal{E}}^1\,\beta  = \beta'\,\ell_{\mathcal{E}}^2, \quad
  \ell_{\mathcal{E}}^2\,\alpha = \alpha'\,\ell_{\mathcal{E}}^1, \qquad
  \beta\,\alpha = \id+s_1\,\ell_{\mathcal{E}}^1, \quad
  \alpha\,\beta = \id+s_2\,\ell_{\mathcal{E}}^2.
\end{equation}
We use the dots $\bullet$ to avoid introducing new notations for the
corresponding spaces of sections of vector bundles.

\begin{definition}
  Two $\mathcal{C}$-differential operators $\Delta_1$ and $\Delta_2$
  on~$\mathcal{E}$ are called \emph{equivalent} if there exist
  $\mathcal{C}$-differential operators $\alpha$\textup{,} $\beta$\textup{,}
  $\alpha'$\textup{,} $\beta'$\textup{,} $s_1$, and $s_2$ such that
  \begin{equation*}
    \Delta_1\,\beta  = \beta'\,\Delta_2, \quad
    \Delta_2\,\alpha = \alpha'\,\Delta_1, \qquad
    \beta\,\alpha = \id+s_1\,\Delta_1, \quad
    \alpha\,\beta = \id+s_2\,\Delta_2.
  \end{equation*}
\end{definition}
(see~\cite{DudnikovSamborski:LOvSPDE} and references therein).
Thus, we can say that the linearization operators $\ell_{\mathcal{E}}^1$ and
$\ell_{\mathcal{E}}^2$ are equivalent.

The following simple Lemma explains why this notion is really important.

\begin{lemma}
  $\mathcal{C}$-differential operators $\Delta_1$ and $\Delta_2$ are 
  equivalent if and only if the $\Delta_1$- and
  $\Delta_2$-coverings are isomorphic as linear coverings.
\end{lemma}

So, to prove that $\ell^*$-covering is invariant we have to establish that the
operators $\ell_{\mathcal{E}}^{1*}$ and $\ell_{\mathcal{E}}^{2*}$ are
equivalent.  This is implied by the following result.

\begin{theorem}
If two normal operators~$\Delta_1$ and~$\Delta_2$ are equivalent then 
$\Delta_1^*$ is equivalent to $\Delta_2^*$.
\end{theorem}

\begin{corollary}
  The equation $\mathcal{L}^*(\mathcal{E})$ does not depend on the embedding
  $\mathcal{E}\to J^{\infty}$.
\end{corollary}

Now, recall that bivectors were defined as conservation laws on
$\mathcal{L}^*(\mathcal{E})$, while operators that correspond to them are
essentially generating functions of these conservation laws.  Thus, the
operators depend on using an embedding $\mathcal{E}\to J^{\infty}$.  Assume
that we have two different embeddings as above, so that they give rise to two
operators $A^1$ and $A^2$ that correspond to the same bivector.  Here are the
formulas that relate these two operators:
\begin{equation}
  \label{a:9}
  \begin{aligned}
    A^2&=\alpha\,A^1\,\alpha'^*, \\
    A^1&=\beta\,A^2\,\beta'^*.
  \end{aligned}
\end{equation}

\section{Examples}

Let us revise the three examples from the Introduction.  

\begin{example}[KdV]
  We considered two different embeddings of the KdV equation to jets:
  \begin{equation*}
    \begin{aligned}
      u_t-u_{xxx}-6uu_x&=0, \\[2ex]
      \begin{pmatrix}
        u_x-v \\
        v_x-w \\
        w_x-u_t+6uv
      \end{pmatrix}&=0.
    \end{aligned}
  \end{equation*}
  Here are all operators of diagram~\eqref{a:3}:
  \begin{equation*}
    \ell_{\mathcal{E}}^1=D_t-D_{xxx}-6uD_x-6u_x,\qquad
    \ell_{\mathcal{E}}^2=
    \begin{pmatrix}
      D_x & -1 & 0 \\
      0 & D_x & -1 \\
      -D_t+6v & 6u & D_x
    \end{pmatrix},
  \end{equation*}

  \begin{equation*}
    \alpha=
    \begin{pmatrix}
      1 \\ D_x \\ D_{xx}
    \end{pmatrix},\quad
    \alpha'=
    \begin{pmatrix}
      \hphantom{-}0 \\ \hphantom{-}0 \\ -1
    \end{pmatrix},\quad
    \begin{aligned}
      \beta&=
      \begin{pmatrix}
        1 & 0 & 0
      \end{pmatrix}, \\
      \beta'&=
      \begin{pmatrix}
        -D_{xx}-6u & -D_x & -1
      \end{pmatrix},
    \end{aligned}
  \end{equation*}
  
  \begin{equation*}
    s_1=0,\quad s_2=
    \begin{pmatrix}
      0 & 0 & 0 \\
      1 & 0 & 0 \\
      D_x & 1 & 0
    \end{pmatrix}.
  \end{equation*}
  Formulas~\eqref{a:9} relate Hamiltonian operators~\eqref{eq:2}
  and~\eqref{a:99}.
  \begin{remark}
    If we take an operator from~\eqref{eq:2} for $A^1$ and compute $A^2$
    via~\eqref{a:9} we
    will get an operator from~\eqref{a:99} only up to the equivalence.
  \end{remark}
\end{example}

\begin{example}[Camassa-Holm equation]
  The Camassa-Holm equation written in the usual form
  $u_t-u_{txx}-uu_{xxx}-2u_xu_{xx}+3uu_x=0$ has a bi-Hamiltonian structure:
  \begin{equation*}
    A_1=D_x\qquad A_2=-D_t-uD_x+u_x.
  \end{equation*}

  If we rewrite the equation in the form
  \begin{align*}
    &m_t+um_x+2u_xm=0, \\
    &m-u+u_{xx}=0
  \end{align*}
  then the bi-Hamiltonian structure takes the form
  \begin{equation*}
    A'_1=
    \begin{pmatrix}
      D_x & 0 \\
      D_x-D_x^3 & 0
    \end{pmatrix}\qquad
    A'_2=
    \begin{pmatrix}
      0 & -1 \\
      2mD_x+m_x & 0
    \end{pmatrix}
  \end{equation*}
  Note that the operators $B_1$ and $B_2$ from Example~\ref{e:1} are entries
  (up to sign) of the matrix $A'_1$ and $A'_2$.  Thus we see that studying
  bi-Hamiltonian structure of the Camassa-Holm equation does not require the
  use of the $(1-D_x^2)^{-1}$ ``operator''.
\end{example}

\begin{example}[Kupershmidt deformation]
  Let $\mathcal{E}$ be a bi-Hamiltonian equation given by $F=0$ and $A_1$ and
  $A_2$ be the Hamiltonian operators.
  \begin{definition}
    The Kupershmidt deformation~$\tilde{\mathcal{E}}$
    of~$\mathcal{E}$ has the form
    \begin{equation*}
      F+A_1^*(w)=0,\qquad A_2^*(w)=0,
    \end{equation*}
    where $w=(w^1,\dots,w^l)$ are new dependent variables.
  \end{definition}
  
  \begin{theorem}
    \label{t:1}
    The Kupershmidt deformation~$\tilde{\mathcal{E}}$ is a bi-Hamiltonian
    system.
  \end{theorem}
  
  The proof of this theorem consists of checking that the following two
  bivectors define a bi-Hamiltonian structure:
  \begin{gather*}
    \tilde{A}_1=
    \begin{pmatrix}
      A_1 & -A_1 \\
      0 & \ell_{F+A_1^*(w)+A_2^*(w)}
    \end{pmatrix}
    \quad \tilde{A}_2=
    \begin{pmatrix}
      A_2 & -A_2 \\
      -\ell_{F+A_1^*(w)+A_2^*(w)} & 0
    \end{pmatrix}
  \end{gather*}

  The generalisation of Kupershmidt's theorem from the Introduction is the
  following.

  \begin{theorem}
    \label{t:2}
    If $\psi_1$,~$\psi_2$,~\dots{} is a Magri hierarchy
    for~$\mathcal{E}$ then, under some technical assumptions,
    $(\psi_i,-\psi_{i+1})$\textup{,}
    $i=1$\textup{,}~$2$\textup{,}~\dots\textup{,} is a Magri hierarchy for the
    Kupershmidt deformation~$\tilde{\mathcal{E}}$.
  \end{theorem}
\end{example}

Details and proofs of Theorem~\ref{t:1} and~\ref{t:2} can be found
in~\cite{KerstenKrasilshchikVerbovetskyVitolo:IKD}.

\begin{acknowledgement}
We wish to thank the organizers and participants of the Abel Symposium~$2008$ 
in Troms\o{} for making the conference a productive and enlivening event. We 
also are grateful to Sergey Igonin for reading a draft of this paper and 
useful comments.
\end{acknowledgement}

\end{document}